\documentclass[12pt]{article}
\usepackage{amsmath,amsfonts,amssymb,theorem}
\usepackage[matrix,arrow]{xy}

 \numberwithin{equation}{subsection}
 \numberwithin{footnote}{section}

 \newtheorem{cor}[subsection]{Corollary}
 \newtheorem{lemma}[subsection]{Lemma}
 
 \newtheorem{thm}[subsection]{Theorem}

 {
 \theorembodyfont{\rmfamily}
 \newtheorem{defn}[subsection]{Definition}
 \newtheorem{example}[subsection]{Example}
 \newtheorem{counterexample}[subsection]{Counter-Example}
\newtheorem{exa-cr}[subsection]{Example--Construction}

 \newtheorem{rem}[subsection]{Remark}

 }
 \newcommand{\qed}{\ifhmode\unskip\nobreak\fi\quad\ensuremath\square}
 \newenvironment{proof}{\paragraph{Proof}}{\par\medskip}

\title{Elements of Nonstandard Algebraic Geometry }

\date{Nov 2001}
\author{Caucher Birkar\thanks{\small School of Mathematical Sciences, Nottingham University, Nottingham, UK.}}

\begin{document}
 \maketitle

\begin{abstract} I investigate the algebraic geometry of nonstandard varieties, using techniques of nonstandard mathematics. 

 \end{abstract}

\tableofcontents
\clearpage


\section{\normalsize{Introduction}}
Methods of nonstandard mathematics have been successfully applied to many
parts of mathematics such as real  analysis, functional analysis, topology, 
probability theory, mathematical physics etc. But just a little bit has been
done in foundations of nonstandard algebraic geometry so far. Robinson indicated some ideas in
[R1] and [R2] to prove Nullstellensatz (R$\ddot{\mbox{u}}$ckert's theorem) and
Oka's theorem, using nonstandard methods, in the case of analytic varieties.

In this paper we try to formulate first elements of nonstandard algebraic geometry. 

Consider an enlargement $^{*}X$ of an affine variety $X$ over an
algebraically closed field $k$. We
often take $k=\mathbb{C}$ to be able to define the shadow of limited points
of $^{*}X$. 

As one of the first results in section 4 (th 4.6) we shall show that the shadow
of any 1-codimensional principal (given by an internal polynomial with a finite
number of monomials) subvariety  of $^{*}X_{^{*}\mathbb{C}}$ is closed in $X$
where  $^{*}X_{^{*}\mathbb{C}}$ is the $^{*}X$ as a variety over the field
$^{*}\mathbb{C}$. 

Also in section 4 (th 4.2) we show that the shadow of any internal open subset
of $^{*}X$ equals $X$, which in turn implies that every point on $X$ has an
internal nonsingular point in its halo.

In section 5 we discuss an error in Robinson's paper [R1, th.5.3]
and indicate a way to fix it.

In section 6 we introduce the notion of a countable
infinite dimensional affine variety and prove Nullstellensatz in the case of
an uncountable underlying algebraically closed field, in particular for the
field of complex numbers.

Finally in section 7 we investigate enlargements of a commutative ring $R$ and 
$R$-modules $M$. We use flatness of $^{*}R$ over $R$ to prove 
$^{*}{M} \simeq {^{*}{R}} \otimes_{R} M
$ for $R$ a Noetherian commutative ring $R$ and a finitely generated
$R$-module $M$.\\

{\bf Acknowledgment} I am grateful to my supervisor I. Fesenko
for his suggestion of the subject  and comments.

\clearpage{\pagestyle{empty}}

\section{\normalsize{List of Notation}}

$ ^{*}\mathbb{C}[z_{1},\dots,z_{n}]$ ....................................
internal polynomials over $ ^{*}\mathbb{C}$ in $n$ variables\\ 
$ (^{*}\mathbb{C})[z_{1},\dots,z_{n}]$ ...............................
polynomials over $^{*}\mathbb{C}$ in $n$ variables\\
  $\mathcal{P}_{F}(A)$
.............................................the set of finite subsets of $A$
\\        $zh_{X}(a)$...............................................Zariski
halo of $a$ in $^{*}X$\\
$h_{X}(a)$.................................................halo of $a$ in
$^{*}X$\\
$^{*}X_{^{*}k}$.................................................$^{*}X$ as a
variety over the field $^{*}k$\\

\clearpage{\pagestyle{empty}}

\section{\normalsize{Basic Definitions}}

  We consider the enlargement of a set which contains an algebraically closed
field $k$ and real numbers. Then we can speak of the enlargement of affine and
projective spaces and more ,the enlargement of any quasiprojective variety.
 Let $X$ be a variety over $k$ and let $^{*}X$ be its enlargement. $^{*}X_{
^{*}k}$ denotes $^{*}X$ as a variety over the field $^{*}k$. Note that this is
completely different from $^{*}X$.
 \begin{defn}
Let $a\in X$ then the halo of $a$ in Zariski topology is 
\[ zh_{X}(a)=\bigcap_{a\in U}{^{*}U}. \] where $U$ is Zariski open in $X$.
\end{defn}
  We distinguish it from $h_{X}
(a)$ which stands for the halo of $a$ when $k=\mathbb{C}$ and $U$ is open
in the sense of usual topology. In this case $^{*}X^{lim}$ shows the
 elements with limited coordinates.\\
 The map $^{*}:X\longrightarrow  {^{*}X}$ is the
natural map which takes $a$ to $^{*}a$ and usually we denote the image of $a$
by the same $a$.
And also we have another important map $^{o}
: {^{*}X}^{lim}\longrightarrow X$ which takes each point to its shadow.\\*
 We would get two different "topologies" on $^{*}X$. One is the internal
Zariski topology such that its opens are the internal open subsets of $^{*}X$.
In order this is not always a topology. That is the intersection of a
collection of closed subsets may not be a closed subset. For example let
$X=\mathbb{A}^{1}_{k}$ and $B_{M}=\{x\in {^{*}\mathbb{N}}: 1\leq x \leq M\}$.
Now let $\mathfrak{B}=\{B_{M}\}_{M \leq N}$ where $N,M$ are unlimited
hypernatural numbers and $k$ is an algebraically closed field with
characteristic $0$. All $B_{M}$ in $\mathfrak{B}$ are hyperfinite and then by
transfer internal closed subsets of $^{*}X$. Now consider $\bigcap_{B\in
\mathfrak{B}}B=\mathbb{N}$ which is not an internal subset of $^{*}X$ and
then not internal closed subset. \\*
   Another topology is the usual Zariski
topology on $^{*}X_{^{*}k}$ as a variety over the field $^{*}k$.

\section{\normalsize{Properties of the $^{*}$ and $^{o}$ maps}}

  $X$ shows an affine variety through this section. 
Consider on $^{*}X$ the internal topology in which a basis of open subsets
consists of complements of all zeros of an internal polynomial
(i.e. an element of $^{*}\mathbb{C}[z]$).

The first thing which draw
our attention is the continuity of the $^{*}$ map. We shall show
that this map is not
continuous.
\begin{example}
Let $X=k=\mathbb{C}$, then there is an internal closed subset of $^{*}X$ with
a non-closed preimage under the $^{*}$ map. 
\end{example}
 The following formula is
true:
 \[(\forall A\in \mathcal{P}_{F}(\mathbb{C}))(\exists p\in
\mathbb{C}[z])(\forall a \in \mathbb {C})(a\in A \longleftrightarrow p(a)=0).
\]  
By transfer we have:

\[(\forall A \in {^{*}\mathcal{P}_{F}(\mathbb{C})})(\exists p\in
 {^{*}\mathbb{C}[z]})(\forall a \in {^{*}\mathbb{C}})(a\in A\longrightarrow
p(a)=0).\]

 Now let $A=\{x\in{^{*}\mathbb{N}}: 1\le x\le N \}$
for an unlimited hypernatural number $N$.
$A$ is  a hyperfinite subset of
$^{*}\mathbb C$.\\
 Then there is an internal polynomial in $^{*}\mathbb{C}[z]$ which
vanishes exactly on $A$. The preimage of $A$ is $\mathbb{N}$ which is not a
closed subset of $\mathbb{C}$.\\
 
We can prove a stronger assertion that for any subset $B$ of $\mathbb{C}$, there is an
internal closed subset of $^{*}\mathbb{C}$ which  has $B$
as its preimage. To
prove it we can consider a hyperfinite approximation of $B$ in
$^{*}\mathbb{C}$, say $H$. $B\subseteq H \subseteq {^{*}B}$. The preimage of
$^{*}B$ is $B$ and then the preimage of $H$ is also $B$.  \\

 Now we look at images of subsets of $^{*}X$ under the $^{o}$ map in the case
of  $k=\mathbb{C}$. Note that we defined the $^{o}$ from $X^{lim}$ to $X$, but
we can consider the image of subsets of $^{*}X$ by taking the images of its
limited points. Unexpectedly the image of any nonempty internal open set is the
whole $X$.\\

\begin{thm}
Let $A$ be a nonempty internal open set in $^{*}X$ then $^{o}A$ is $X$.
\end{thm}
\begin{proof}
It is sufficient to prove the theorem for principal internal open subsets.
Then let $A= {^{*}X}_{f}$ be a nonempty internal principal open subset where
$f$ is an internal polynomial. If the shadow of $A$ is not $X$ there should be
some point $a \in X$ for which $f(h_{X}(a))=0$. Hence the following formula is
true:

\[(\exists g\in {^{*}\mathbb{C}[z_{1},\dots,z_{n}]})(\exists \varepsilon \in
{^{*}\mathbb{R}^{+}})(\forall z\in {^{*}X})((\exists w \in {^{*}X})(g(w)\neq 0
 ) \wedge (\mid z-a \mid \leq \varepsilon\rightarrow g(z)=0)).\]
And by transfer:
\[(\exists g\in \mathbb{C}[z_{1},\dots,z_{n}])(\exists \varepsilon \in
\mathbb{R}^{+})(\forall z\in X)((\exists w \in X)(g(w)\neq 0) \wedge (\mid
z-a\mid \leq \varepsilon\rightarrow g(z)=0)).\]
 It is easy to see that the latter is 
not true. \end{proof}

\begin{cor}
There is a nonsingular point $\xi$ in $h_{X}(a)$ for every $a\in X$. 
\end{cor}

\begin{thm}
Let $f:X\longrightarrow Y$ be a regular map of varieties over
$\mathbb {C}$. Then we have:
\begin{enumerate}
\item (i)  $^{o}(^{*}Z)=Z$ for every closed subset of $X$;
\item (ii)  $(^{*}f)^{-1}(^{*}Z)= {^{*}(f^{-1}(Z)})$ for every  
subset $Z$ of $Y$. 
\end{enumerate}
\end{thm}
\begin{proof}
(i) Obviously $Z\subseteq {^{o}(^{*}Z)}$. 
Let $Z =  V (g_{1}, \dots ,g_{l} )$ and
 let $x\in {^{o}(^{*}Z)}$,
$x={^{o}\xi}$ with $\xi \in {^{*}Z}$.
Then $g_{i}(\xi)=0$ for $1\leq i \leq l$. Clearly  $g_{i}(^{o}\xi)=0$ for
$1\leq i \leq l$ and this proves that $x\in Z$.\\
(ii) consider the formula :
\[(\forall x\in X)(x\in (f)^{-1}(Z)\longleftrightarrow f(x)\in Z).\]
and by transfer:
\[(\forall x\in {^{*}X})(x\in {^{*}((f)^{-1}(Z)})\longleftrightarrow
{^{*}f(x)}\in {^{*}Z}).\]

and on the other hand we have:
 \[(\forall x\in {^{*}X})(x\in(^{*}f)^{-1}(^{*}Z)\longleftrightarrow
{^{*}f(x)}\in {^{*}Z}).\]   which proves the equality.\\
\end{proof}
 It is well known that the shadow of any subset of $^{*}\mathbb{R}$ is a
closed subset in $\mathbb{R}$, the field of real numbers, in the sense of
real topology. But that is not such easy in the case of algebraic sets. Now we
show that the shadow of an internal closed subset of $^{*}X$ is not always
closed in $X$. For example consider $B_{M}$ in
$^{*}\mathbb{A}^{1}_{\mathbb{C}}$, which was introduced in section 1, with $M$
an unlimited hypernatural number. Obviously $^{o}B_{M}=\mathbb{N}$ which is not
closed in $\mathbb{A}^{1}_{\mathbb{C}}$. A better deal is to consider closed
subsets of $^{*}X_{\mathbb{C}}$. 

\begin{thm}
Let $f\in (^{*}\mathbb{C})[z_{1},\dots,z_{n}]$ be a polynomial with limited
coefficients and let $^{o}f$ be nonzero. Then we have:
\[
 ^{\circ}(V (f) ) ={V}( ^{\circ}f).
\]
\end{thm}
\begin{proof}
The shadow of $f$, $^{\circ}f$, may be a constant i.e the coefficients of
nonzero degree monomials in $f$ are infinitesimal. This implies that no
limited point could be in $V (f)$. On the other hand ${V}( ^{\circ}f
)=\emptyset$. Then the equality is proved in this case.
 Otherwise let ${\xi} \in {
^{*}\mathbb{C}^{n} }$ be a limited point and $ f(\large\xi )=0$. Then $
^{o}f(^{o} \large \xi)= 0 $ hence $^{o}\large\xi\in V( ^{o}f ) $.\\  Now let $a\in
V(^{o}f)$ then $f(a)\simeq 0$. Using hypotheses,
$f(h_{\mathbb{C}^n}(a))\subseteq h_{\mathbb{C}^n}(0)$. It is sufficient to
find a point in the halo of $a$ such that $f$ vanishes at that point. Now if
$f(a)\neq 0$ we can change variables linearly such that $a$ is transferred to
origin. Note that the new polynomial, say $g$ has also limited coefficients
and this translation takes $h_{\mathbb{C}^{n}}(a)$ to $h_{\mathbb{C}^{n}}(0)$.
We have $g= g^{inf}+g^{ap}$ such that  $g^{inf}$ has infinitesimal
coefficients and $g^{ap}$ has appreciable coefficients. Then
$^{o}g={^{o}g^{ap}}$.\\* Now we use induction on the number of variables. If
$n=1$ Robinson--Callot theorem[DD, ch. 2,th. 2.1.1] shows that
$g(h_{\mathbb{C}}(0))=h_{\mathbb{C}}(0)$ because $g$ is S-continuous as it has 
limited coefficients. If $1<n$ we consider the form with highest degree
appeared in $g^{ap}$, say $h$. $h$ is a sum of monomials of the same degree.\\* 
If $h=\alpha z_{1}\dots z_{n}$ where $\alpha$ is a hypercomplex number, then we
change variables such that $z_{1}=w_{1}$ and $z_{i}=w_{i}+w_{1}$. This change,
obviously maps the halo of origin on itself and we get a new polynomial $e$ with  
limited coefficients from $g$. Now consider $e(w_{1},\dots,w_{n-1},0)$,
clearly the shadow of this polynomial in a smaller than $n$ number of
variables, is not constant, and we use induction.\\* 
In the remaining cases we can again replace  one of the variables by zero and
reduce the number of variables, if necessary, and use induction.\\ This proves the existence
of a zero for $f$ and completes proof of the theorem.  \end{proof}

 We can generalize this result by replacing $\mathbb{C}^{n}$ with its
affine subvariety, $X$. The theorem is again true. Although the previous
theorem is a particular case of the next theorem, but their proofs are
of different nature and we prefer to keep the previous proof. 
\begin{thm}
Let $X$ be an algebraic closed subset of $\mathbb{C}^{n}$ and $f\in
(^{*}\mathbb{C})[z_{1},\dots,z_{n}]$ be a polynomial with limited
coefficients. Then there is a $g$ in $(^{*}\mathbb{C})[z_{1},\dots,z_{n}]$
with limited coefficients such that:  \[ V (g)= V(f),  
  {^{\circ}V(f)} =V(^{\circ}g).
\]
where zeros of these polynomials are taken in $^{*}X$ and $X$
correspondingly.
 \end{thm}
\begin{rem}
It is not always true if we take $g$ to be $f$ itself. For example let
$X=V(z_{1})$ in $\mathbb{C}^{2}$ and $f=z_{1}+\varepsilon z_{2}$ in which
$\varepsilon$ is an infinitesimal hypernatural number. Then $^{o}f=z_{1}$
which is identically zero on $X$. But the shadow of $V(f)$ is just a
single point.
\end{rem}
\begin{proof}
If $V(f)={^{*}X}$ then the theorem is trivial.\\*
 In other cases 
if $^{\circ}f$ is not identically zero on $X$ we take $g=f$. Otherwise let
 $\ddot{f}$ be $f$ divided by one of its  coefficients with maximum absolute
value. If $^{o}\ddot{f}(X)\neq0$ then put $g=\ddot{f}$. Otherwise assume that
$^{o}\ddot{f}(^{*}X)=0$ then $V(\ddot{f}-{^{o}\ddot{f}})=V(\ddot{f})=V(f)$.
$\ddot{f}-{^{o}\ddot{f}}$ has smaller number of monomials  than $f$. By
continuing this process eventually we get a polynomial $g$ such that its
shadow is not identically zero on $X$ and $V(g)=V(f)$.\\*
 Now let $x\in {^{\circ}V(g)}$, then
$x={^{o}\xi}$ for some $\xi \in V (g)$. From $g(\xi)=0$ we get 
$^{o}g(^{o}\xi)=0$ and then $x\in V (^{o}g)$.  
 Conversely let $x\in V(^{\circ}g)$. We want to prove that
$h_{X}(x)\cap V(g)\neq \emptyset$. Let $Y\subseteq X$ be an irreducible curve
containing $x$. It is sufficient to prove that $ h_{Y}(x)\cap V(g)\neq
\emptyset$. It is proved if $g(^{*}Y)=0$, otherwise change variables such that
$x$ be transferred to the origin and then consider  $^{*}Y_{^{*}\mathbb{C}}$.
$V(g)\cap {^{*}Y_{^{*}\mathbb{C}}}$ is a finite set i.e a zero dimensional
subvariety, say $\mathcal{A}=\{\xi_{1},\dots,\xi_{l}\}$. Since $^{*}X\subseteq
{^{*}\mathbb{C}^{n}}$, then every point of $^{*}X$ is as
$(b_{1},\dots,b_{n})$, with $n$ coordinates, $b_{1},\dots,b_{n}$. Now if no
point in $\mathcal{A}$ is infinitesimal, with infinitesimal coordinates,  then
every $\xi_{i}$ has at least a non-infinitesimal coordinate,
say $a_{i_{j}}$. The index $j$ means that $a_{i_{j}}$ has appeared in the
$j$-th coordinate of $\xi_{i}$. Now put
$h_{i}=(z_{j}-a_{i_{j}})/a_{i_{j}}$. And  let $h=h_{1}\times \dots \times
h_{l}$. Now obviously $\mathcal{A}\subseteq V(h)$. Then we have $h^{t}=eg$ on
$^{*}Y$, for some polynomial $e$ and natural number $t$. By construction $h$
and $g$ have limited coefficients. $e$  should also have limited coefficients,
otherwise $h^{t}/s=(e/s)g$ on $^{*}Y$ where $s$ is a coefficient appeared in
$e$ with maximum absolute value. Then $^{o}(h^{t}/s)=0={^{o}(e/s)}{^{o}g}$ on
$Y$. But $Y$ is irreducible, hence $^{o}(e/s)=0$ on $Y$. Now we can use the
method by which we constructed $g$ and reduce the number of monomials appeared
in $e$. Then we get a new $e$ with limited coefficients which satisfies 
$^{o}e\neq 0$, $h^{t}=eg$ and  $^{o}h^{t}={^{o}e}{^{o}g}$. This is a
contradiction, because $^{o}h$ is not zero at origin.      
\end{proof}

\section{\normalsize{Generic Points for Prime Ideals}}

 Let $\Gamma$ be the ring of analytic functions at origin (origin of 
$\mathbb{C}^{n}$). An important theorem in complex analysis says
that every prime ideal of $\Gamma$ has a generic point in the halo of
origin.
 We prove a similar theorem in the algebraic context.
\begin{thm}
Let $X$ be an irreducible affine variety and $x \in X$. Then every prime ideal
in the ring of regular functions at $x$ has a generic point in the Zariski
halo of $x$. 
\end{thm}
\begin{proof}
 
Let $\mathfrak{p}$ be a prime ideal in $\mathcal{O}_{X,x}$, the ring of
regular functions at $x$. Define:
 \[A_{f,g,U}=\{y\in U: U \mbox{ is open in $X$}
\wedge g(y)\neq 0 \wedge f(y)=0 \wedge \mbox{$f,g$ are regular on $U$}\}.\]

 Using Nullstellensatz $A_{f,g,U}\neq\emptyset$ where $x\in U$, $f\in
\mathfrak{p}$ and   $g\notin \mathfrak{p}$. Similarly the collection
$\{A_{f,g,U}\}_{x\in U,f\in\mathfrak{p},g\notin \mathfrak{p}}$ has finite
intersection property. Then there would be a $\xi$ in the following set: 
\[\bigcap_{x\in
U,f\in\mathfrak{p},g\notin \mathfrak{p}} {^{*} A_{f,g,U}}.\]
 So $\xi$ is a generic
point for $\mathfrak{p}$ and $\xi\in zh_{X}(x)$ . 
\end{proof}

Thus, we deduce  that the map $\pi : zh_{X}(x) \longrightarrow
Spec(\mathcal{O}_{X,x})$ is surjective where $\pi(\xi)=m_{\xi}$, elements
of $\mathcal{O}_{X,x}$ vanishing at $\xi$. This map  demonstrates how close
$zh_{X}(x)$ and $Spec(\mathcal{O}_{X,x})$ are.\\

\begin{thm}
With the hypotheses of the previous theorem we get:
\[\pi^{-1}(V_{S}(I))= V_{zh}(I).\]
in which $I$ is an ideal of $\mathcal{O}_{X,x}$ , $V_{S}(I)$ is the
closed subset of  $Spec(\mathcal{O}_{X,x})$ defined by $I$ and $
V_{zh}(I)$ is the zeros of $I$ in $zh_{X}(x)$.  
\end{thm}
\begin{proof}
Let $\xi\in zh_{X}(x)$ and $\pi (\xi)\in V_{S}(I)$. Then obviously
$I\subseteq \pi(\xi)$, in other words  every member of $I$ vanishes at $
\xi$. This shows that $\xi$ is in the right side of the above equality.\\*
Conversely let $\xi$ be in the right side of the equality then every member
of $I$ vanishes at $\xi$. This implies that $I$ is contained in $\pi(\xi)$
i.e $\xi$ is in the left side of the equality. 
\end{proof}
 In the analytic case the existence of the generic point is used to prove the
Nullstellensatz theorem. That is if $f,g_{1},\dots,g_{l} \in \Gamma$ and 
${\large V}(g_{1},\dots,g_{l})\subseteq V(f)$ then some power of $f$ should be in
the ideal generated by $g_{i}$'s [R1, sect. 4]. In [R1, th.5.1] the existence
of a generic point was proved for infinite dimensional spaces
$\mathbb{C}^{\Lambda}$, in which $\Lambda$ is an
arbitrary infinite set. Robinson used the previous result
to deduce Nullstellensatz in
this case [R1, th.5.3]. Unfortunately, his prove is erroneous. Now we
indicate the gap.\\

$\large \mbox{Analysis of Robinson's Proof.}$ Let $\Gamma$ be the set of
cylindrical analytic functions in the origin of  $\mathbb{C}^{\Lambda}$  
each one depending only on a finite number of variables. Let
$\mathcal{A}\subseteq \Gamma$ be such that $V(\mathcal{A})\subseteq 
{\large V}(f)$ in a neighborhood of origin. 
If no power of $f$ is in $<\mathcal{A}>$
then there is a prime ideal, say ${P}$ containing $\mathcal{A}$ and
not $f$. ${P}$ has a generic point in the halo of origin, say $\xi$.
Robinson  concludes that $f$ is zero at $\xi$ because  $
{\large V}(\mathcal{A})\subseteq V(f)$ in a neighborhood of origin like $U$. But
this is not true. Consider:  

\[(\forall x\in U)((\forall h\in \mathcal{A}) h(x)=0\longrightarrow
f(x)=0).\]   and by transfer:   
                                                                              
\[(\forall x\in {^{*}U})((\forall h\in {^{*}\mathcal{A}}) h(x)=0\longrightarrow
{^{*}f(x)}=0).\] 
This formula is true but it is different from:

\[(\forall x\in {^{*}U})((\forall h\in {^{im}\mathcal{A}})h(x)=0\longrightarrow
{^{*}f(x)}=0).\] 
 which is a wrong formula  Robinson applied to $\xi$.
\begin{counterexample}
Let $\Lambda=\mathbb{C}$, $h_{a}=z_{a}(z_{0}-a)-1$, $\mathcal{A}=\{h_{a}:
a\in \mathbb{C}\mbox{ \it and } a\neq 0\}$ and $f=z_{0}$ in which
$z_{a}$ is a variable indexed by  $a$. Then $V(\mathcal{A})\subseteq 
{\large V}(f)$ and no power of $f$ is in $<\mathcal{A}>$.
\end{counterexample}
Let $\xi \in V(\mathcal{A})$, then $z_{0}(\xi)=0$ because for every
nonzero $a\in\mathbb{C}$, $h_{a}(\xi)=z_{a}(\xi)(z_{0}(\xi)-a)-1=0$ and then
$z_{0}(\xi)-a$ is nonzero. Hence $z_{a}(\xi)=1/(z_{0}(\xi)-a)=1/(-a)$. This means
that $V(\mathcal{A})=\{\xi\}$. Clearly $ \xi\in V(z_{0})$. But if a
power of $z_{0}$, say $z_{0}^{l}$, be in $<\mathcal{A}>$ then
$z_{0}^{l}=\sum_{i=1}^{t}e_{i}h_{a_{i}}$
 where $h_{a_{i}}\in \mathcal{A}$. Now we can find a point
at which all $h_{a_{i}}$'s are zero and $z_{0}$ is not. But this is a
contradiction. Then no power of $z_{0}$ is in $<\mathcal{A}>$.

\section{\normalsize{Varieties of Infinite Dimensions}}

The previous section demonstrates some peculiar features of  
 varieties of infinite
dimensions. In this section at first we show that Nullstellensatz does not hold
in infinite dimensional algebraic geometry as well as in
infinite dimensional complex analysis.

\begin{counterexample}There is a set $\Lambda$ and a proper ideal $\mathfrak{J}$ in
$S$, the ring of polynomials over $\mathbb{C}$ in variables indexed by
$\Lambda$, such that $V(\mathfrak{J})=\emptyset$.
  \end{counterexample}
Let $\Lambda=\mathbb{C}\cup \{\mathbb{C}\}$, $h_{a}=z_{a}(z_{0}-a)-1$ for
$a\neq 0$ in $\mathbb{C}$ and  $h_{\mathbb{C}}=z_{\mathbb{C}}z_{0}-1$.
Let $\mathfrak{J}$ be the ideal generated by all these functions in $S$.
 Then ${\large V}(\mathfrak{J})=\emptyset$. If $\mathfrak{J}=S$ then there are $a_{1},\dots,a_{l}$ ($a_{l}$ can
be $\mathbb{C}$) and $f_{1},\dots,f_{l}$ such that
\[\sum_{i=1}^{l}f_{i}h_{a_{i}}=1.\]  Now consider all variables which occur in
this formula and let $R$ be the ring of polynomials  in
these variables over $\mathbb{C}$ and $\mathbb{C}^{m}$ the corresponding affine space. Then the
ideal generated by $h_{a_{1}},\dots, h_{a_{l}}$ in $R$, is $R$ itself. That is
$V(\mathfrak{J})=\emptyset$ in $\mathbb{C}^{m}$. This  is not
possible, because we can find a point in  $\mathbb{C}^{m}$ at  which all
$h_{a_{i}}$'s are zero. But right side of the above equation would not be
zero at that point.\\

Fortunately this is not the end of the story. We prove a complete version of
Nullstellensatz similar to the finite dimensions, in the particular case of   $\Lambda=\mathbb{N}$. Let $S$ be the ring
$\mathbb{C}[z_{1},z_{2},\dots]$.
 
\begin{defn} 
Let $X \subseteq \mathbb{C}^{\mathbb{N}}$. We say $X$ is an affine variety in
 $\mathbb{C}^{\mathbb{N}}$ if $X= V(\mathfrak{J})$ for some ideal
$\mathfrak{J}$ of $S$ and we call $\mathbb{C}[X]=S/I(X)$  the ring of regular
functions on $X$. Similarly the field of fractions of $\mathbb{C}[X]$ denoted
by $\mathbb{C}(X)$ is called the field of rational functions on $X$. 
\end{defn} 

\begin{thm}
Let $\mathfrak{M}$ be a maximal ideal of $S$. Then $
{\large V}(\mathfrak{M})\neq\emptyset$.
\end{thm}
\begin{proof}
If for every $n\in\mathbb{N}$ there be a  $a_{n}\in\mathbb{C}$ such that
$z_{n}-a_{n}\in\mathfrak{M}$ then $ \mathfrak{M}=<z_{n}-a_{n}>_{n\in
\mathbb{N}}$, because $<z_{n}-a_{n}>_{n\in\mathbb{N}}$ is a maximal ideal of
$S$. Hence $V(\mathfrak{M})=\{(a_{n})_{n\in\mathbb{N}}\}$.
 Now suppose there is a $n\in\mathbb{N}$ such that 
$z_{n}-a\notin\mathfrak{M}$ for any  $a\in\mathbb{C}$. For simplicity we can
take $n=1$. Now let $S_{i}=\mathbb{C}[z_{1},\dots,z_{i}]$ and
$\mathfrak{M}_{i}$ the contraction  of $\mathfrak{M}$ in  $S_{i}$. 
$\mathfrak{M}_{i}$ is a prime ideal in $S_{i}$ but our goal is to prove that it
is also a maximal ideal.\\*
Let $Y_{i}=V(\mathfrak{M}_{i})$ in $\mathbb{C}^{i}$. Then by our hypothesis
$Y_{1}=\mathbb{C}$,  i.e $\mathfrak{M}_{1}=0$. For every $i$ we have a
projection:
 \[\pi_{i}: Y_{i}\longrightarrow \mathbb{C}.\]
Where $\pi_{i}(y_{1},\dots,y_{i})=y_{1}$.
 Every member of $S$ is a polynomial with a finite number of variables
occurred in it. Then $\bigcup \mathfrak{M}_{i}= \mathfrak{M}$. By a theorem in
algebraic geometry [SH, ch. I,\S 5,th.6] $\pi_{i}( Y_{i})$ is open in
$\mathbb{C}$ or a point in it. If  $\pi_{i} (Y_{i})$ is just a point for some
$i$, say $b$, then $z_{1}-b \in \mathfrak{M}_{i}$ which is a contradiction.
 If all  $\pi_{i} (Y_{i})$ are open, let $x \in\mathbb{C}$. Then there is a
$h\in S$ such that $1-h(z_{1}-x)\in \mathfrak{M}$ and then $1-h(z_{1}-x)\in
\mathfrak{M}_{j}$ for some $j$. $x$ can't be in $\pi_{j}(Y_{j})$ because
$1-h(z_{1}-x)$ doesn't vanish at any point where its coordinate corresponding
to $1$ is $x$. This proves the following equality: \[\mathbb{C}=
\bigcup_{i=1}^{\infty} {\mathbb{C}\setminus\pi_{i}}(Y_{i}).\] which  is
impossible.    \end{proof}
 
This theorem shows that every proper ideal of $S$ at least has a zero in
$\mathbb{C}^{\mathbb{N}}$.
\begin{cor}
An ideal $\mathfrak{M}$ in $S$ is maximal iff it is as $<z_{i}-a_{i}>_{i\in
\mathbb{N}}$ for some $a_{i}\in \mathbb{C}$.
\end{cor} 
 In the proof of the previous theorem we haven't used any specific property of
$\mathbb{C}$, we have just used that it is algebraically closed and 
uncountable. So  
 \begin{cor} 
 The theorem holds if we replace
$\mathbb{C}$ by any uncountable algebraically closed field $k$.
\end{cor}

Now we look at other parts of Nullstellensatz. 
 \begin{thm}
Let $\mathfrak{J}$ be an ideal in $S$, then $I(V(\mathfrak{J}))=\surd
\mathfrak{J}$.
\end{thm}
\begin{proof}
One inclusion  is obvious. Put $T=\mathbb{C}^{\mathbb{N}}$ and let $
{\large V}(\mathfrak{J})\subseteq V(g)$ where $g \in S$. Now we consider a new
space of the same shape, say $W=
\mathbb{C}\times T$. 
We will have a new variable like $z_{0}$ and a new
coordinate corresponding to it (note that $0\notin \mathbb{N}$ in this work). Consider the
ideal $\mathfrak{J^{+}}=\mathfrak{J}+<1-z_{0}g>$ in the ring $S[z_{0}]$.
$\mathfrak{J^{+}}$  has no zero in $W$, so $\mathfrak{J^{+}}=S[z_{0}]$. Hence
there are $h_{0},h_{1},\dots, h_{l}$ in $S[z_{0}]$ and $f_{1},\dots,f_{l}$ in
$\mathfrak{J}$ for which we have: 
\[\sum_{i=1}^{l}h_{i}f_{i}+h_{0}(1-z_{0}g)=1.\]
Now we can put $z_{0}=1/g$ and conclude that either $\mathfrak{J}=S$ or
some power of $g$ is in $\mathfrak{J}$.
\end{proof}
\begin{cor}
Let $\mathfrak{J_{1}},\mathfrak{J_{2}}$ be ideals in $S$ then we have the
following:
\begin{description}
\item[(i)] $V(\mathfrak{J_{1}} \mathfrak{J_{2}})=V(\mathfrak{J_{1}} \cap
\mathfrak{J_{2}})=V(\mathfrak{J_{1}}) \cup V(\mathfrak{J_{2}})$;\\
\item[(ii)] $V(\mathfrak{J_{1}}+ \mathfrak{J_{2}})=V(\mathfrak{J_{1}}) \cap
V(\mathfrak{J_{2}})$;\\
\item[(iii)] $ \surd \mathfrak{J_{1}}$ is prime iff  $V(\mathfrak{J_{1}})$ is
irreducible.
\end{description}
\end{cor}
\begin{proof}
Standard. 
\end{proof}
\begin{defn}
Let $\phi : X \longrightarrow Y$ be a map in which $X$ and $Y$ are
affine varieties. $\phi$ is a regular map if $\phi=(\phi_{1},\phi_{2},\dots)$
in which $\phi_{i}$ is a regular function on $X$.
 Similarly if all $\phi_{i}$ are rationals on $X$ and $\phi(Dom(\phi))
\subseteq Y$, $\phi$ is called a rational map. 
\end{defn}

It is easy to check that, there is an equivalence between the category of
affine varieties over $\mathbb{C}$ (as defined in 4.2) and the category of
reduced countably generated $\mathbb{C}$-algebras.

 It is not obvious that every rational map has a nonempty domain. 
\begin{thm} 
$Dom(\phi)\neq \emptyset$ for any rational map $\phi:X \longrightarrow Y$.
\end{thm}
\begin{proof}
Let $\phi=(\phi_{1},\phi_{2},\dots)$, $\phi_{i}=g_{i}/f_{i}$ and
$T=\mathbb{C}^{\mathbb{N}}$. It is sufficient to prove that there is a point
at which none of $f_{i}$'s vanishes. Suppose there is no such point
i.e. 
\[\bigcup _{i=1}^{\infty} V(f_{i})= \mathbb{C}^{\mathbb{N}}.\]
 Now let $W= \mathbb{C} \times \mathbb{C} \times \mathbb{C} \times \mathbb{C}
\dots$. We define a coordinate system on $W$ such that the
$(2i-1)$th component in it is same as the $i$th component of $T$ i.e we
associate the variable $z_{i}$ to the component with number $2i-1$, and the
variable $w_{i}$ to the $2i$th component.\\*

 Now consider the set:
\[\mathcal{A}=\{1-w_{i}f_{i}: i\in \mathbb{N}\}.\]
 This set has no zero in $W$. Then by theorem 5.3, 
$<\mathcal{A}>=\mathbb{C}[z_{1},w_{1},z_{2},w_{2},\dots]$ and hence there are
$h_{1},\dots,h_{l}$ in $\mathbb{C}[z_{1},w_{1},z_{2},w_{2},\dots]$
such that:
\[\sum_{j=1}^{l}h_{j}(1-w_{i_{j}}f_{i_{j}})=1.\]
 But this is a contradiction because we know that there is some $\xi \in T$
such that $f_{i_{j}}(\xi)\neq0$ for $1\leq j \leq l$. By putting
$w_{i_{j}}(\xi)=1/f_{i_{j}}(\xi)$ we get a point at which all
$(1-w_{i_{j}}f_{i_{j}})$ are zero.
\end{proof}
\begin{cor}
Neither $\mathbb{C}^{\mathbb{N}}$ or $\mathbb{C}^{n}$ ($n$ is finite) is 
the union of a countable set of proper subvarieties.
\end{cor}
 We just proved it for $\mathbb{C}^{\mathbb{N}}$. Suppose that
$\mathbb{C}^{n}=\bigcup _{i=1}^{\infty} V(f_{i})$
in which $f_{i}$ is in $\mathbb{C}[z_{1},\dots,z_{n}]$. Now extend it to
$\mathbb{C}^{\mathbb{N}}$ and we get the result.\\

 Let $S_{\mathbb{N}}=\mathbb{C}[z_{1},z_{2},\dots]$ and
$S_{i}=\mathbb{C}[z_{1},\dots,z_{i}]$. We have inclusions when $n<m$:
\[S_{n} \longrightarrow S_{\mathbb{N}}\]
\[S_{n} \longrightarrow S_{m}\]
and by transfer we have 
\[^{*}S_{n} \longrightarrow {^{*}S_{\mathbb{N}}}\]
\[^{*}S_{n} \longrightarrow S_{N}\]
in which $S_{N}=\mathbb{C}[z_{1},\dots,z_{N}]$ is the set of internal
polynomials over $^{*}\mathbb{C}$ in variables $z_{1},\dots,z_{N}$ with  an
unlimited hypernatural number $N$.\\*
 Now let $\mathfrak{J}$ be an ideal in $S_{\mathbb{N}}$, $\mathfrak{J}_{n}$
its contraction in $S_{n}$ and $\mathfrak{J}_{N}$ the corresponding
internal ideal in $S_{N}$. We have the following diagram:

$$  \xymatrix{
    S_{n} \ar[d]\ar[r]^{\alpha_{n,\mathbb{N}}} 
 & S_{\mathbb{N}} \ar[d] \\
    {^{*}S_{n}} \ar[r]^{\alpha_{n,N}} 
 &  S_{N}\ar[d] \\
 & {^{*}S_{\mathbb{N}}} \\
}$$ \\

 By using transfer we can see that
$\alpha_{n,N}^{-1}(\mathfrak{J}_{N})={^{*}\mathfrak{J}_{n}}$, for all $n\in
\mathbb{N}$. And then
$\alpha_{\mathbb{N},N}^{-1}(\mathfrak{J}_{N})=\mathfrak{J}$.\\*  

\section{\normalsize{Enlargement of Commutative Rings}}

 In this section we study the enlargement of commutative rings, especially Noetherian rings. 
In the theory of commutative rings localization and completion of rings and modules have 
some typical properties like preserving exactness of sequences and their
closed relation with  tensor product.  That is, if $R$ is Noetherian ring,
$\mathfrak{p}$ a prime ideal  and $M$ is a finitely generated $R$-module then
we have:  
 \[M_{\mathfrak{p}}\simeq R_{\mathfrak{p}} \otimes_{R} M .\]
 \[\widehat{M} \simeq \widehat{R} \otimes_{R} M .\]
     
 We prove similar properties of enlargement of modules. As usual we denote the enlargement 
 of $R$ and $M$ as $^{*}R$ and $^{*}M$. For any ideal $\mathfrak{I}$ of $R$
we have   two notions of radical of $^{*}\mathfrak{I}$ in the ring $^{*}R$.
One is the usual   $\surd ^{*}\mathfrak{I}$ when we consider $^{*}R$ as a ring
and another is the   internal notion of radical, say $^{int}\surd
^{*}\mathfrak{I}$ which is exactly the enlargement of   $\surd \mathfrak{I}$
i.e.

 \[^{int}\surd ^{*}\mathfrak{I}={^{*}\surd \mathfrak{I}}. \]
 Similarly we have the same situation for many other notions. 
 From now on we work with  a Noetherian commutative ring  $R$.

\begin{thm}
For any ideal $\mathfrak{I}$ in $R$ we have: 
\[^{*}min(\mathfrak{I})=min_{int}(^{*}\mathfrak{I})=min(^{*}\mathfrak{I}).\]
\end{thm}
\begin{proof}
$R$ is Noetherian then $min(\mathfrak{I})$ is a finite set, say
$\{\mathfrak{p}_{1},\dots,\mathfrak{p}_{l}\}$.  Then it is its own enlargement. Now let $\mathfrak{q}$ be a prime ideal in $^{*}R$   containing
$^{*}\mathfrak{I}$.   Hence its contraction  $\mathfrak{q}^{c}$ in $R$ is a
prime ideal containing $\mathfrak{I}$.    There is some $j$ such that
$\mathfrak{p_{j}} \subseteq \mathfrak{q^{c}}$. Then     $^{*}\mathfrak{p_{j}}
\subseteq \mathfrak{q}$. This implies the equalities.    
\end{proof}
 It can also easily be proved that $^{*}J(R)=J(^{*}R)$ where $J(R)$ is the
Jacobson radical of $R$ and similarly $J(^{*}R)$ is the Jacobson radical of
$^{*}R$.  \begin{cor}.\\
(i) $^{int}\surd ^{*}\mathfrak{I}=\surd ^{*}\mathfrak{I}$ and
$nil_{int}(^{*}R)={^{*}nil(R)}= nil(^{*}R)$;\\
(ii) $\mathfrak{q}$ is $\mathfrak{p}$-primary iff $^{*}\mathfrak{q}$ is
$^{*}\mathfrak{p}$-primary  iff $^{*}\mathfrak{q}$ is
internally  $^{*}\mathfrak{p}$-primary. 
\end{cor}
Let $\phi : M \longrightarrow N$ be a
homeomorphism of $R$-modules. Then  
\begin{lemma}
 (i) $ker {^{*}\phi}= {^{*}ker \phi}$;\\
 (ii) $im {^{*}\phi}= {^{*}im} \phi$.\\
 \end{lemma}
 \begin{proof}
 (i) \[(\forall m\in M) ( m\in \ker \phi \longleftrightarrow \phi(m)=0).\]
     and by transfer:
     \[(\forall m\in {^{*}M}) ( m\in {^{*}\ker \phi} \longleftrightarrow
{^{*}\phi}(m)=0).\] (ii) Use a similar formula.
\end{proof} 
 \begin{cor}
 Let $M,N,L$ and $K$ be $R$-modules. Then \\*
 (i) $0 \longrightarrow N \longrightarrow M \longrightarrow K \longrightarrow 0$
 is exact iff 
  $0 \longrightarrow {^{*}N} \longrightarrow {^{*}M} \longrightarrow {^{*}K}
\longrightarrow 0$ is exact;\\
(ii) $^{*}M/^{*}K={^{*}(M/K)}$;\\
\end{cor}

\begin{lemma}
$^{*}R$ is a faithfully flat $R$-algebra. 
\end{lemma}
\begin{proof}
 By [B,ch. I,\S 2,$n^{o}$11]  $^{*}R$ is a faithfully flat $R$-algebra iff for
any maximal ideal $\mathfrak{m}$   in $^{*}R$, $\mathfrak{m} ^{*}R \neq
{^{*}R}$ and any solution of an $R$-homogeneous linear equation 
$\sum_{i=1}^{l}a_{i}Y_{i}=0$ in $^{*}R^{l}$ is an $^{*}R$ linear combination
of solutions in $R^{l}$. \\*
  Let $\mathfrak{m}$ be  any maximal ideal of $R$. 
 Since $R$ is Noetherian, $ \mathfrak{m} ^{*}R={^{*}
\mathfrak{m}}$, then 
 $\mathfrak{m} ^{*}R \neq {^{*}R}$.   \\* 
  Now let $f=\sum_{i=1}^{l}a_{i}Y_{i}=0$ be an $R$-homogeneous
linear equation. Let $\mathcal{A}$    be the module of solutions to $f$ in
$R^{l}$. $\mathcal{A}$ is an $R$-submodule of $R^{l}$.    Since $R$ is
Noetherian then $\mathcal{A}$ is finitely generated, say
$\mathcal{A}=<\beta_{1},\dots,\beta_{c}>$.  Then we have:

  \[(\forall x_{1},\dots,x_{l}\in R)[\sum_{i=1}^{l}a_{i}x_{i}=0
\longleftrightarrow (\exists r_{1},\dots,r_{c} \in
R)(x_{1},\dots,x_{l})=\sum_{i=1}^{c}r_{i}\beta_{i} ]. \] 
and using transfer:  
 \[(\forall x_{1},\dots,x_{l}\in {^{*}R})[\sum_{i=1}^{l}a_{i}x_{i}=0
\longleftrightarrow (\exists r_{1},\dots,r_{c} \in
{^{*}R})(x_{1},\dots,x_{l})=\sum_{i=1}^{c}r_{i}\beta_{i} ].\]  This proves that
$^{*}R$ is $R$-flat, and then faithfully flat $R$-algebra.  \end{proof}  

Let $M$ be a finitely generated $R$-module.
Define a bilinear function   
\[\omega : M \times {^{*}R} \longrightarrow
{^{*}M}\]
 such that $\omega (m,r)=rm$. This induces a unique
$R$-homomorphism 
\[ \Omega_{M} : M \otimes_{R} {^{*}R} \longrightarrow {^{*}M}, \quad     
\Omega_{M}(\sum_{i=1}^{t}a_{i}(m_{i} \otimes
r_{i}))=\sum_{i=1}^{t}a_{i}r_{i}m_{i}.\]
 where $a_{i}\in R$, $m_{i}\in M$
and $r_{i}\in {^{*}R}$. Clearly $\Omega$ is surjective.   
 \begin{thm}
 $\Omega_{M}$ is an isomorphism.
 \end{thm}
 \begin{proof}
  We first assume that $M$ is a free module, say $M=R^{s}$. Let
$\{e_{1},\dots,e_{s}\}$    be a basis for $M$ over $R$. Then every element of
$M \otimes_{R} {^{*}R}$ can be written as    $\sum_{i=1}^{s}a_{i}(e_{i} \otimes
r_{i})$ and its image under $\Omega_{M}$ will be   
$\sum_{i=1}^{s}a_{i}r_{i}e_{i}$. Now assume $\sum_{i=1}^{s}a_{i}r_{i}e_{i}=0$.
     By transfer all $a_{i}r_{i}$ should be zero. This proves the theorem
when $M$      is free.\\*      Now in the general case, there is an $l$ and a
surjective homomorphism       from $R^{l}$ to $M$. Let $K$ be the kernel of
this homomorphism. Then we get an       exact sequence of $R$-modules:
     \[0 \longrightarrow K\longrightarrow R^{l} \longrightarrow M
\longrightarrow 0.\]     and so
      \[0 \longrightarrow {^{*}K} \longrightarrow {^{*}R^{l}} \longrightarrow
{^{*}M} \longrightarrow 0.\]        And also by flatness of $^{*}R$ we have:
      \[0 \longrightarrow K \otimes_{R} {^{*}R}\longrightarrow R^{l}
\otimes_{R} {^{*}R} \longrightarrow M \otimes_{R} {^{*}R} \longrightarrow
0.\]    Now the maps $\Omega$, namely $\Omega_{K}$, $\Omega_{R^{l}}$ and
$\Omega_{M}$    give us vertical homomorphisms between the two exact
sequences:    

$$  \xymatrix{
 0 \ar[r] 
 & K \otimes_{R} {^{*}R} \ar[d]\ar[r]^{\lambda} 
 & R^{l} \otimes_{R} {^{*}R} \ar[d]\ar[r]^{\gamma} 
 & M \otimes_{R} {^{*}R} \ar[d]\ar[r] 
 & 0 \\
0 \ar[r] 
 & {^{*}K} \ar[r]^{\alpha} 
 & {^{*}R^{l}} \ar[r]^{\beta} 
 & {^{*}M} \ar[r] 
 & 0\\
}$$ \\
 
Suppose $\Omega_{M}(a)=0$. There is $b$ such that $\gamma(b)=a$. And let $\Omega_{R^{l}}(b)=c$. 
By commutativity of the diagram $\beta(c)=0$. Hence there is $d$ such that
$\alpha(d)=c$. $\Omega_{K}$ is surjective then there is $e$ such that
$\Omega_{K}(e)=d$. Then $\Omega_{R^{l}}\lambda(e)=c$. But $\Omega_{R^{l}}$ is
an isomorphism, then $\lambda(e)=b$. And by exactness
$\gamma(b)=\gamma\lambda(e)=0$. This shows that $\Omega_{M}$ is an isomorphism
of $R$-modules. This completes the proof. 
\end{proof}

 We can consider $M \otimes_{R} {^{*}R}$ as a $^*R$-module. $\Omega_{M}$ is
also $^*R$-homomorphism and then it is a $^*R$-isomorphism.\\
 
By [B, ch.
IV,\S2.6,th.2] we have  \[Ass_{^{*}R}{^{*}M}=Ass_{^{*}R}(M \otimes_{R}
{^{*}R})=\{^{*}\mathfrak{p}: \mathfrak{p} \in Ass_{R}M \}.\] 

\begin{cor}
$Ass_{^{*}R}{^{*}M}={^{*}Ass_{R}M}$.
\end{cor}

 By [VS,th1.1] we can say that (as a particular case)
$T={^{*}\mathbb{C}}[z_{1},\dots,z_{m}]$ is a faithfully flat  $
S=(^{*}\mathbb{C})[z_{1},\dots,z_{m}]$-algebra. By lemma $6.6$,
$T$ is also a faithfully flat $R=\mathbb{C}[z_{1},\dots,z_{m}]$-algebra:

$$  \xymatrix{
 S \ar[r]^{\beta} 
 & T \\
 R \ar[u]^{\gamma} \ar[ru]^{\alpha}
 & \\
}$$ \\

 Now let
$\mathfrak{J}$ be an ideal of $R$. $(\mathfrak{J}S)T={^*\mathfrak{J}}$ and 
let $\mathfrak{J}_{1}=\gamma^{-1}(\mathfrak{J}S)$.  By flatness of $\beta$,
$\beta^{-1}(^{*}\mathfrak{J})= \mathfrak{J}S$, hence
$\gamma^{-1}[\beta^{-1}(^{*}\mathfrak{J})]=\mathfrak{J}_{1}$. On the other
hand by flatness of $\alpha$, $\alpha^{-1}(^{*}\mathfrak{J})=
\mathfrak{J}$, then we conclude that
$\gamma^{-1}(\mathfrak{J}S)=\mathfrak{J}$. Then we get another diagram:

$$  \xymatrix{
 S/\mathfrak{J}S \ar[r]^{\beta} 
 & T/ ^{*}\mathfrak{J}\\
 R/\mathfrak{J} \ar[u]^{\gamma} \ar[ru]^{\alpha}
 & \\
}$$ \\

\begin{cor}
$\mathfrak{J}$ is prime iff $\mathfrak{J}S$ is prime.
\end{cor}

  If $\mathfrak{J}$ be a radical ideal. Then $\mathfrak{J}S$ and
$^{*}\mathfrak{J}$ are also radical.\\*
 These ideals correspondingly define closed subsets $Y$ (in
$\mathbb{C}^{m}$), $^{*}Y_{^{*}\mathbb{C}}$ (in
$^{*}\mathbb{C}^{m}_{^{*}\mathbb{C}}$) and $^{*}Y$ (in $^{*}\mathbb{C}^{m}$).
Moreover $R/\surd\mathfrak{J}, S/\surd\mathfrak{J}S$ and $T/
\surd^{*}\mathfrak{J}$ are their coordinate rings. Now using the previous
corollary we get 

\begin{cor}
$^{*}Y$ is irreducible iff $^{*}Y_{^{*}\mathbb{C}}$ is irreducible
iff $^{*}Y$ is internally irreducible. \end{cor}

\section*{\normalsize{References:}}
\begin{description}
\item[AM] M.F.Atiyah, I.G.Macdonald; Introduction to Commutative Algebra,
Addison-Wesley,1969.\\
\item[B] N.Bourbaki, Elements of Mathematics, Commutative Algebra. Herman 1972.\\ 
\item[DD] F.Diener, M.Diener; Nonstandard Analysis in Practice, Springer-Verlag 
1995.\\
\item[R1] A.Robinson; Germs, Applications of Model Theory to Algebra, Analysis,
and Probability (ed.W.A.J.Luxemburg) New York, etc. 1969 pp. 138-149.\\ 
\item[R2] A.Robinson; Enlarged Sheaves, Lecture Notes in Mathematics 369, pp.
249-260 Springer-Verlag 1974.\\ 
\item[R3] A.Robinson;
Nonstandard Analysis, North-Holland,1974.\\
\item[Sh] I.Shafarevich; Basic Algebraic Geomerty. 
Springer-Verlag,1972.\\ 
\item[VS] L.Van den Dries, K.Schmidt; Bounds in the theory of
polynomial rings over fields, A nonstandard approach.Inventiones mathematicae,
Springer-Verlag, 1984.

\end{description}

\end{document}